\documentclass[12pt, twoside, draft]{article}
\usepackage[cp1251]{inputenc}
\usepackage{enumerate}

\usepackage{ifthen}
\usepackage{amsmath}
\usepackage{amsfonts}
\usepackage{latexsym}
\usepackage{amsthm}
\usepackage{amssymb}

\newcounter{chapter}
\newcommand{\logo}{\baselineskip2pc \hbox to\hsize{\hfil\copyright\,\footnotesize
\CopyName, \Year}}

\theoremstyle{plain}
\newtheorem{theorem}{Theorem}
\newtheorem*{MTheorem}{Main Theorem}

\newtheorem*{corollary}{Corollary}
\newtheorem*{theomKL}{Theorem K--L}

\theoremstyle{definition}

\theoremstyle{remark}
\newtheorem{remark}{Remark}

\newcommand{\sep}{/\kern-2pt/ }

\newcommand{\vs}{\vspace{.1in}}
\newcommand{\vsk}{\vspace{.2in}}

\newcounter{lang}
\newcommand{\received}[1]
{\hfill{\footnotesize \textit{%
  \ifthenelse{\value{lang}=0}{Received}
    
#1}}}

\newcommand\matrefrus{{\ /\kern-2pt/ ???. ???\-?.\ -- \rm\Year.\ -- ?.\Volume, \No \Number.\
-- C.\Pages.}}

\newcommand{\CopyName}{B.\ N.\ Khabibullin} 
\newcommand{\NAME}{B.\ N.\ KHABIBULLIN} %
\newcommand{\Year}{2010} 
\newcommand{\Volume}{34} 
\newcommand{\Number}{2} 
\newcommand{\Pages}{197--206} 
\newcommand{\rightheadtext}{NEVANLINNA'S THEOREMS} 
\usepackage[english,russian,ukrainian]{babel}

\renewcommand{\refname}{\refnam}
\newtheorem*{theomN}{Nevanlinna's Theorem}

\DeclareMathOperator{\Zero}{Zero}
\DeclareMathOperator{\Exh}{Exh}
\DeclareMathOperator{\Hol}{Hol}
\DeclareMathOperator{\Harm}{harm}
\DeclareMathOperator{\SH}{sbh}
\DeclareMathOperator{\Mer}{Mer}
\DeclareMathOperator{\Av}{\mathrm A}

\newcommand{\R}{\mathbb{R}}
\newcommand{\bC}{\mathbb{C}}
\newcommand{\N}{\mathbb{N}}

\newcommand{\D}{\mathbb{D}}
\newcommand{\e}{\varepsilon}
\DeclareMathOperator{\Int}{int}
\DeclareMathOperator{\dist}{dist}
\DeclareMathOperator{\supp}{supp}

\DeclareMathOperator{\Pol}{Pol}
\renewcommand{\d}{{\rm \, d}}

\vs
\newcommand{\tit}{NEVANLINNA'S THEOREMS} 
\date{}
\begin{document}
\hbox to \textwidth{\footnotesize\textsc  Математичні Студії. Т.\Volume, \No \Number
\hfill
Matematychni Studii. V.\Volume, No.\Number}
\vspace{0.3in}
\markboth{{\NAME}}{{\rightheadtext}}
\begin{center} \textsc {\CopyName} \end{center}
\begin{center} \renewcommand{\baselinestretch}{1.3}\bf {\tit} \end{center}

 {{\bf Abstract.}{\small { 
 B.\ N.\ Khabibullin,
\textit{Nevanlinna's theorems}

We give a short survey on generalizations of Nevanlinna's theorems on zero distribution of bounded holomorphic functions and representation  
of meromorphic  functions 
in multiply connected
domains. It is a part of our report     
in the conference on complex analysis dedicated to the  memory of Anatolii Asirovich Goldberg in Lviv, May 31-June 5, 2010.
}} \vsk

\renewcommand{\refname}{\refnam}

\vskip10pt

{\bf Definitions, agreements, and basic notions}. 
Let $\Omega$ be a domain in the complex plane $\mathbb C$ or in the Riemann sphere $\mathbb C_{\infty}\neq \Omega$.
For $S\subset \Omega \subset \bC_{\infty}$, we denote by $\overline S$ and  $\partial \Omega$ 
the {\it closure\/} and the {\it boundary\/}  
of $S$ relative to $\bC_\infty$. 
We write $S\Subset \Omega$ if $\overline S \subset \Omega$.

Denote  by $\Hol(\Omega)$, $\Mer (\Omega)$, $\SH (\Omega)$, and $\Harm (\Omega)$  the classes of all holomorphic, mero\-m\-o\-r\-ph\-ic, subharmonic, and harmonic functions on $\Omega$.

We are concerned with finite or infinite sequences $\Lambda =\{ \lambda_k \}$, $k=1,2, \dots$ of not necessarily distinct points from the domain $\Omega$, without limit points in $\Omega$. 
Let $n_{\Lambda}$ be an integer-valued {\it counting measure\/} of sequence $\Lambda$ defined by 
\begin{equation*}
n_{\Lambda}(S):=\sum_{\lambda_k \in S} 1	, \quad S\subset \Omega .
\end{equation*}

Let $S\subset \Omega$. $\Lambda \subset S \; \Longleftrightarrow\; \supp n_{\Lambda}\subset S$.

The sequence $\Lambda$ coincides with a sequence $\Gamma =\{\gamma_n\}$
(or is equal to $\Gamma$, or $\Lambda =\Gamma$)
iff $n_{\Lambda}=n_{\Gamma}$. $\Lambda \subset \Gamma$
means $n_{\Lambda}\leq n_{\Gamma}$.
$\Lambda \cap \Gamma$ and $\Lambda \cup \Gamma$ are defined by 
$n_{\Lambda \cap \Gamma}:=\min \{ n_\Lambda , n_\Gamma\}$ and
$n_{\Lambda \cup \Gamma}:=n_\Lambda + n_\Gamma$.

Given $f\colon A\to B$ and $b\in B$ , we write {\it $f\equiv b$ on $A'$} if $f$ is identically equal to $b$ on $A' \subset A$; in the opposite case,  $f\not\equiv b$ on $A'$. 

Let $A,B\subset [- \infty, +\infty ]$.
A function $f\colon A\to B$ is {\it increasing\/} ({\it decreasing} resp.) if, for any $x_1, x_2\in A$,  $x_1\leq x_2$ implies  $f(x_1)\leq f(x_2)$
\; ( $f(x_1)\geq f(x_2)$ resp.). 

Given $a\in \R$, and $f\colon A\to [-\infty , +\infty]$, we set $a^+:=\max\{0,a\}$, $f^+:=\max \{0, f \}$.
  
The term ``positive'' (``negative'' resp.) means ``$\geq 0$'' (``$\leq 0$'' resp.).

Let $f\in \Hol(\Omega)$ or $f\in \Mer (\Omega)$, $f\not\equiv 0, \infty$ on $\Omega$. Write  $\Zero_f $ for the zero set  of $f$ (counting multiplicities). Evidently, $\Zero_f$ is 
a sequence of not necessarily distinct points from the domain $\Omega$, without limit points in $\Omega$. 

A sequence $\Lambda$ is  the {\bf zero sequence} for a subspace $H\subset \Hol(\Omega)$
(further we write $\Lambda \in \Zero(H)$) if and only if there exists a function $f\in H$ such that $\Lambda = \Zero_f$.

A function $f\in \Hol(\Omega)$ {\it vanish on\/} $\Lambda$ if and only if  $\Lambda \subset \Zero_f$ (we write $f(\Lambda)=0$).

A sequence $\Lambda$ is the  {\bf zero \underline{sub}sequence} or the {\bf non-uniqueness sequence} for the space $H$  iff there exists a nonzero function $f\in H$  such that  $f(\Lambda)=0$.

{\bf Problems}. Let $H$ be a subspace of $\Hol(\Omega)$. We consider the following five  problems:
\begin{enumerate}[\rm 1.] {\it 
\item What point sequences $\Lambda$ can be {\bf zero sequences  for} $H$?

\item What point sequences $\Lambda$ can be {\bf zero \underline{sub}sequences  for} $H$?

\item In what cases is a  zero 
{\bf \underline{sub}sequence  for} $H$ simultaneously {\bf zero sequence for} $H$ or for some, preferably minimal, {\bf extension space} $\widehat H \supset H$ of holomorphic functions in $\Omega$?
\item When can\/ {\bf a meromorphic function} in $\Omega$ be represented as a ratio of\/ {\bf holomorphic functions from\/} $H$?
\item When can\/ {\bf a meromorphic function} in $\Omega$ be represented as a ratio of  holomorphic functions from\/ $H$ {\bf without common zeros}?      
}
\end{enumerate}

{\bf Green's function}. Denote by 
$$
g_{_\Omega} (\cdot , z)\colon \bC_{\infty}\setminus \{z\}\to [0,+\infty)
$$ 
the {\it extended Green's function\/} for $\Omega$ with a pole at $z\in \Omega$, i.\,e., $g_{_\Omega}(\zeta , z )\equiv 0$ for points $\zeta \in \bC_{\infty} \setminus \overline \Omega$, 
$g_{_\Omega}(\cdot , z)\in \SH ( \bC_{\infty} \setminus \{ z \})$,  
$g_{_\Omega}(\cdot , z)\in \Harm ( \Omega \setminus \{ z \})$, and also 
$$
g_{_\Omega}(\zeta , z)=-\log |\zeta -z |+O(1), \quad \zeta \to z.
$$ 

Given a continuous function $\phi \colon \partial \Omega \to \mathbb R$, we denote by $H_{_\Omega}\phi$  the {\it solution the Dirichlet problem\/} for $\Omega$ with boundary function $\phi$ or the associated {\it Perron function\/} 
\begin{equation*}
H_{_\Omega}\phi:=\sup\{u \in \SH (\Omega)\colon \\
\limsup_{z\to \zeta} u(z)\leq \phi (\zeta), \; \forall \zeta\in \partial \Omega  \}.
\end{equation*}

{\bf Harmonic measure}. Denote by $\mathcal B (\partial \Omega)$  the $\sigma$-algebra of Borel subset of $\partial \Omega$. 

Denote by $\omega_{_\Omega}(z, \cdot)$ the {\it harmonic measure\/} 
for $\Omega$ at the point $z\in \Omega$, i.\,e., 
$$
\omega_{_\Omega}(\cdot , \cdot) \colon \Omega\times \mathcal B (\partial \Omega) \to [0,1]
$$  
 such that
\begin{enumerate}[{a)}]
\item the map $B\mapsto \omega_{_\Omega} (z, B )$ is a Borel probability measure  on $\partial \Omega$;
 \item if $\phi \colon \partial \Omega \to \mathbb R$ is continuous function, then 
 $$
 H_{_\Omega}\phi (z)=\int_{\partial \Omega} \phi (\zeta) \d\omega_{_\Omega}(z, \zeta).
 $$
\end{enumerate}

{\bf The unit disk}. Starting points of our researches are Nevalinna's theorems (1929).

Denote by $\Hol^{\infty}(\D ) \subset \Hol (\D )$  the space of holomorphic bounded functions on $\D$.

Denote by $f_{\Lambda}$ a holomorphic function $\not\equiv  0$ on $\D$ with zero sequence $\Zero_{f_{\Lambda}}=\Lambda\subset \Omega$.

\begin{theomN}[{\rm on zeros}]  The following 
 statements are equivalent:
\begin{enumerate}[{\rm 1)}]
	\item $\Lambda$ is a zero sequence for\/   $\Hol^{\infty}(\D)$;
	
	\item $\Lambda $ is a zero \underline{sub}sequence for\/  $\Hol^{\infty}(\D)$;
	
\item for \underline{each} or \underline{some} function\/ $f_{\Lambda} \in \Hol (\D)$ 
\begin{multline*}
\sup_{r<1}\bigl(H_{_{r\D}}\log|f_{\Lambda}|\bigr)(0)=\sup_{r<1}\int_{r\partial\D}\log \bigl|f_{\Lambda}(z)\bigr| d\omega_{_{r\D}}(0, z)\\=\sup_{r<1} \frac1{2\pi} \int_0^{2\pi} \log \bigl|f_{\Lambda}(re^{i\theta})\bigr| \d \theta <+\infty \, ;	
\end{multline*}

\item $\sum\limits_k g_{_{\D}}(\lambda_k, 0)= \sum\limits_k\log\dfrac1{|\lambda_k|}<+\infty$
$\Longleftrightarrow \; 
\sum\limits_k (1-{|\lambda_k|})<+\infty $.

\end{enumerate}
\end{theomN}

Let $z\in \Omega$, $p,q \in \Hol (\Omega)$, and let
\begin{equation}\label{r:f}
f=\frac{p}{q}\, \in \Mer(\Omega), \;  p(z)=q(z)=1. 
\end{equation}
 We set
\begin{equation}\label{f:N}
u_f:=\max \{ \log |p|, \log |q|\}\in \SH (\Omega) \; 
\end{equation}
\begin{equation}\label{f:AS}
\text{ or }
u_f
:=\log \sqrt{ |p|^2 +|q|^2}\in \SH (\Omega).
\end{equation}

Let $D$ be a subdomain of $\Omega$, $z\in D \Subset \Omega$. The integral
\begin{equation}\label{f:T}
T_{_D}(f; z):=\int_{\partial D} u_f \d \omega_{_D}(z, \cdot)
\end{equation}
is {\it Nevanlinna's characteristic\/} of $f$ on $D$ (in  the form of Ahlfors--Shimizu  for \eqref{f:AS}) at the point $z$ relative to the representation \eqref{r:f}.

If $\Pol_f=\{\gamma_k\}_{k=1}^{\infty}\subset \Omega$ is the {\it pole sequence\/} of $f$ in $\Omega$, i.\,e., zero sequence of $1/f$ in $\Omega$, and $\Zero_p\cap \Zero_q=\varnothing$, then, for \eqref{f:N},  
\begin{equation*}
T_{_D}(f;z)=\sum_{k}g_{_D}(\gamma_k, z)\\+\int_{\partial D} \log^+|f| \d \omega_{_D}(z,\cdot ).
\end{equation*}
If $z=0\in D$, then $T_{_D}(f):=T_{_D}(f;0)$. 

{\bf Multiply connected domains}. There are generalizations of Nevanlinna's theorems to  classes of holomorphic and meromorphic functions on special finitely connected domains $\Omega$.

Given $z\in \bC$ and $0<t<+\infty$, we write
 $$D(z,t):=\{ w\in \bC \colon |w-z|<t\}, \quad \overline D (z,t):= \overline{D(z,t)}.$$ 

We consider now results from [1] (Andriy~Kondratyuk, and Ilpo~Laine, 2006).
The authors develop the Nevalinna theory and connected topics for meromorphic and holomorphic functions to $(m+1)$-connected domains $\Omega$, $m\in \N$, for following cases.

\begin{enumerate} 

\item An $(m+1)$-connected domain 
$
\Omega=\bC \setminus \bigcup_{j=1}^{m}\{c_j\}\quad, c_j\in \bC,  \quad m\in \N,
$ is called the {\it $m$-punctured plane.\/} For example, such 2-connected  domain is the {\it punctured plane\/} $\bC \setminus \{0\}$.

 \item A bounded\footnote{always bounded relative to $\bC$} $(m+1)$-connected domain
   $$
   \Omega= D(0,R)\setminus \Bigl(\bigcup_{j=1}^m D(z_j, r_j)\Bigr), \; m\in \N,
   \quad D(z_j, t_j)\Subset  D(0,R),$$ $\overline D(z_j, t_j)\cap \overline D(z_{j'}, t_{j'})=\varnothing$ for all  $j\neq j'$, is called {\it the strictly circular domain\/}.  {2-connected} annuli 
  $
  A_R=D(0,R)\setminus D(0,1/R), \quad 1<R<+\infty,
  $ 
  are examples of such domains.
\item Let $\Omega \subset \bC$ be a bounded $(m+1)$-connected domain, $m\in \N$, 
$$
\Omega = \bigcup_{j=0}^{m} G_j,
$$ 
where $G_j$ are simply connected domains (relative to $\bC_{\infty}$),  $\partial G_j$ are simple (Jordan) closed paths, $G_0$ is bounded, and $\infty \in G_j$ for $j=1, \dots , m$,  $\bC\setminus G_j\Subset G_0$, 
$$
(\bC\setminus G_j)\bigcap (\bC\setminus G_{j'})=\varnothing ,
\; j\neq j', \; j,j'\geq 1.
$$
Such domains $\Omega$ are called in [1] {\it \underline{admissible}}.\/
Let $\varphi_j$ are conformal mappings 
of $\D$ onto $G_j$ realizing a conformal equivalence of $\D$ and $G_j$ with $\varphi_0(0)=0$, and $\varphi_j (0)=\infty$ as $j=1, \dots, m$. 
Then every $\varphi_j$ can be extended to a homeomorphism of $\overline \D$ onto $\overline G_j$ (the Carath{\'e}odory Theorem). For a simple closed path $\Gamma$, $\Int \Gamma$ denote the bounded domain with the boundary $\Gamma$ (the Jordan Theorem).

Let 
$$
\Gamma_{jr}(s):=\varphi_{j}(e^{is}), \; 0\leq s\leq 2\pi ,\; j=0,\dots , m  ,
\quad \Gamma_{jr}^*:=\Gamma_{jr}\bigl( [0, 2\pi)\bigr), \; j=0,\dots , m ,
$$
and
$$
\Omega_r:=\bigl(\Int \Gamma_{0r}^*\bigr)\setminus \bigcup_{j=1}^k \overline{\Int \Gamma_{jr}^*}, \; r_0\leq r<1.  
$$
Let  $\Lambda=\{ \lambda_k\}$ be a sequence in $\Omega$, and let
$$
N_0 (r, \Lambda):=\int_{r_0}^r
\frac{n_\Lambda(\Omega_t)}{t} \d t,
\; r_0\leq r<1,
$$
where $r_0<1$ is a constant sufficiently close to $1$.

For a meromorphic function $f$ on an \underline{admissible} bounded domain $\Omega$ denote 
\begin{equation*}
m_0 (r,f):=\\
\frac1{2\pi} \sum_{j=1}^m \biggr(\int_0^{2\pi}\log^+ \bigl|f\bigl(\varphi_j(re^{is})\bigr)\bigr|\d s\biggr.\\
-\biggl.\int_0^{2\pi}\log^+ \bigl|f\bigl(\varphi_j(r_0 e^{is})\bigr)\bigr|\d s\biggr), 
\end{equation*}
$r_0\leq r<1.$ The function 
$$
T_0(r,f):=N_0(r, \Pol_f)+m_0(r,f), \; r_0\leq r<1,
$$
is the {\it Nevalinna characteristic\/} of $f$ (in the sense of A.~Kondratyuk and I.~Laine). 
\end{enumerate} 

\begin{theomKL}[{\rm on meromorphic functions, [1; Theorem 43.2]}] Let $f\in \Mer (\Omega)$, where $\Omega\subset \bC$ is a \underline{finitely connected}  \underline{bounded} \underline{admissible}  domain. 
If 
$
T_0(r,f)=O(1),\quad  r\to 1,
$ 
then there are \underline{bounded} functions $h_1, h_2\in \Hol (\Omega)$ such that $f=h_1/h_2$.   
\end{theomKL}

 Let $f_\Lambda\in \Hol (\Omega)$ be a 
function with zero sequence $\Lambda \subset \Omega$. If we apply this theorem to the function $f_\Lambda\in \Hol(\Omega)$, then we get
\begin{theomKL}[{\rm on zeros}] Let $\Omega\subset \bC$ be a \underline{finitely connected}  \underline{admissible} \underline{bounded}  domain. If 
$$
\sup_{r_0\leq  r<1}m_0(r,f_\Lambda)<+\infty,
$$ 
then 
there is a \underline{bounded} holomorphic on $\Omega$ function $f$ such that $f(\Lambda)=0$.  		

\end{theomKL}

We consider the Problems 1--5 in a general form in [2], [3].

A function (weight) $M \colon \Omega \to \R$  define weighted space 
$\Hol(\Omega; M){:=}$
\begin{equation*}
 \left\{ f\in \Hol(\Omega) \colon 
\sup_{z\in \Omega} \frac{|f(z)|}{\exp M(z)}<+\infty
\right\} \, .
\end{equation*}
If $M\equiv 0$ on $\Omega$, then $\Hol (\Omega;0)=\Hol^{\infty}(\Omega)$.

Let $\Omega$ be a domain in $\bC_{\infty}$.  
 
We use following classification of domains $\Omega$ such that
$
\boxed{0\in \Omega \not\ni \infty}.
$
\begin{enumerate}[{\bf I.}]
\item This domain $\Omega$ is 
\begin{enumerate}[{a)}]
\item {\bf simply connected} (for example, $\Omega =\mathbb C$ or $\Omega=\mathbb D$) or   
\item {\bf  finitely connected}  and $\boxed{\overline \Omega \neq \mathbb C_{\infty}}$  (for example, each \underline{bounded} finitely connected domain $\Omega$).
\end{enumerate}
\item This domain $\Omega$ is {\bf $(m+1)$-connected} with $m\in \mathbb N$ and $\boxed{\overline\Omega = \mathbb C_{\infty}}$.

\item This domain $\Omega$ is
\begin{enumerate}[{a)}]
	\item {\bf bounded} or
	\item {\bf unbounded}.
	
\end{enumerate}

\end{enumerate}

Let $\Exh \Omega=\{D\}$ be a {\it exhaustion\/} of $\Omega$, i.\,e., $\cup_{D \in \Exh \Omega} D=\Omega$, where $0\in D$ and every $D\in \Exh \Omega$ is the regular domain for the Dirichlet problem. Such exhaustion exists always (countable increasing and such that all $D$ have smooth boundary).

\begin{theorem}[{\rm on zeros}] Let $\Lambda \subset\Omega$, and let $\Omega$ be a domain of type\/ {\bf I}. 
Following statements are equivalent:
\begin{enumerate}[{\rm 1)}]
	\item $\Lambda$ is a zero sequence for\/   $\Hol^{\infty}(\Omega)$;
	
	\item $\Lambda $ is a zero \underline{sub}sequence for\/  $\Hol^{\infty}(\Omega)$;
	
\item for \underline{each} or \underline{some} function\/ $f_{\Lambda} \in \Hol (\Omega)$ with\/ $\Zero_{f_{\Lambda}}=\Lambda$ 
\begin{equation*} 
\sup_{D\in \Exh \Omega}\bigl(H_{_{D}}\log|f_\Lambda|\bigr)(0):=\sup_{D\in \Exh \Omega}\int_{D}\log \bigl|f_{\Lambda}(z)\bigr| d\omega_{_D}(0, z)<+\infty \, ;	
\end{equation*}

\item $\sup\limits_{D\in \Exh \Omega} \sum\limits_k g_{_{D}}(\lambda_k, 0)<+\infty$.
\end{enumerate}
If\/ $\Omega$ possesses the Green's function, then we can remove  $\boxed{\sup\limits_{D\in \Exh \Omega}}$ everywhere and replace $D$ with $\Omega$.
\end{theorem}

\setcounter{theorem}{0}
\begin{theorem}[{on meromorphic functions}] Let  $\Omega$ be a domain of type\/ {\bf I}.
Let\/ {\rm (see \eqref{r:f})} $$
f =\frac{p}{q}\in \Mer (\Omega), \; p,q \in \Hol (\Omega), \; p(0)=q(0)=1,$$
and\/ {\rm (see  \eqref{f:T}, \eqref{f:N}, \eqref{f:AS} resp.)}
\begin{equation*}
T_{_D}(f):=\int_{\partial D} u_f \d \omega_{_D}(0, \cdot), \; D\in \Exh(\Omega),
\end{equation*}
where
\begin{equation*}
u_f:=\max \{ \log |p|, \log |q|\}\in \SH (\Omega) \; 
\end{equation*}
\begin{equation*}
\text{ or }
u_f
:=\log \sqrt{ |p|^2 +|q|^2}\in \SH (\Omega).
\end{equation*}

Assume that one of the following two conditions holds:
\begin{enumerate}[{\rm 1)}]
\item $\sup\limits_{D\in \Exh \Omega} T_{_D}(f) < +\infty$;
\item $p,q\in \Hol^{\infty}(\Omega)$.
\end{enumerate}
Then there are  $p_0,q_0\in \Hol^{\infty}(\Omega) $ {\bf without common zeros} such that $f=p_0/q_0$.
\begin{remark} If\/ the domain $\Omega$ is regular for the Dirichlet problem, then we can remove  $\boxed{\sup\limits_{D\in \Exh \Omega}}$ in\/ {\rm 1)} and replace $D$ with $G$.
\end{remark}
\end{theorem}

\begin{theorem}[{\rm on zeros}] Let $\Lambda \subset\Omega$, and let $\Omega$ be a domain of type\/ {\bf II}. If  
	$\Lambda$ is a zero sequence for\/   $\Hol^{\infty}(\Omega)$,
then
\begin{enumerate}	
	\item[{\rm 2)}] $\Lambda $ is a zero \underline{sub}sequence for\/  $\Hol^{\infty}(\Omega)$;
	
\item[{\rm 3)}] for \underline{each} (for \underline{some}) function $f_{\Lambda} \in \Hol (\Omega)$ with\/ $\Zero_{f_{\Lambda}}=\Lambda$ 
\begin{equation*} 
\sup_{D\in \Exh \Omega}\bigl(H_{_{D}}\log|f_\Lambda|\bigr)(0):=\sup_{D\in \Exh \Omega}\int_{D}\log \bigl|f_{\Lambda}(z)\bigr| d\omega_{_D}(0, z)<+\infty \, ;	
\end{equation*}

\item[{\rm 4)}] $\sup\limits_{D\in \Exh \Omega} \sum\limits_k g_{_{D}}(\lambda_k, 0)<+\infty$.
\end{enumerate}

Conversely, assume that one of the conditions {\rm 2)--4)} holds.
Then there is a constant $b<m$ such that $\Lambda$ is zero sequence for every space $\Hol (\Omega ,M )$
with
\begin{equation}\label{df:Ml}
M\colon z\mapsto c_0^+\log^+|z|+
\sum_{k=1}^m c_k^+\log^+\frac1{|z-a_k|}\, ,  
\end{equation}
where $\sum_{k=0}^mc_k=b$, $a_k\in \bC \setminus \Omega$.
\end{theorem}

\setcounter{theorem}{1}
\begin{theorem}[{on meromorphic functions}] Let  $\Omega$ be a domain of type\/ {\bf II}.
Let\/ {\rm (see \eqref{r:f})} $$
f =\frac{p}{q}\in \Mer (\Omega), \; p,q \in \Hol (\Omega), \; p(0)=q(0)=1.$$

Assume that one of the following two conditions holds:
\begin{enumerate}[{\rm 1)}]
\item $\sup\limits_{D\in \Exh \Omega} T_{_D}(f) < +\infty$;
\item $p,q\in \Hol^{\infty}(\Omega)$.
\end{enumerate}
Then there are constant $b<m$ such that, for every function $M$ from \eqref{df:Ml}, there exist functions  
$p_0,q_0\in \Hol(\Omega; M)$ {\bf without common zeros} representing $f=p_0/q_0$.

\end{theorem}

\begin{theorem}[{\rm on zeros}] Let $\Omega\subset \bC$ be a domain, and let $\Lambda \subset\Omega$ be a sequence. If  
	$\Lambda$ is a zero (sub)sequence for\/   $\Hol^{\infty}(\Omega)$,
then
\begin{enumerate}	

\item[{\rm 3)}] for \underline{each} or \underline{some} function\/ $f_{\Lambda} \in \Hol (\Omega)$ with\/ $\Zero_{f_{\Lambda}}=\Lambda$ 
\begin{equation*} 
\sup_{D\in \Exh \Omega}\bigl(H_{_D}\log|f_\Lambda|\bigr)(0):=\sup_{D\in \Exh \Omega}\int_{D}\log \bigl|f_{\Lambda}(z)\bigr| d\omega_{_D}(0, z)<+\infty \, ;	
\end{equation*}
		
\item[{\rm 4)}] $\sup\limits_{D\in \Exh \Omega} \sum\limits_k g_{_{D}}(\lambda_k, 0)<+\infty$.
\end{enumerate}

Conversely, let  one of the conditions\/ {\rm 3), 4)} is fulfilled.
Then $\Lambda$ is zero \underline{sub}sequence for  space $\Hol (\Omega ,M )$
with 
$$
M \colon z\mapsto \log \frac1{\dist (z, \partial \Omega)}\,,\quad z\in \Omega,
$$
where $\dist(z, \partial \Omega)$ is the Euclidean distance from $z$ up to $\partial \Omega$, if $\Omega$ is \underline{bounded},\\ and with any
\begin{equation*}
M \colon z\mapsto \log \frac1{\dist (z, \partial \Omega)}\\+c_0\log^+|z|+c_1\log^+\frac1{|z-a|}\,,\quad z\in \Omega,
\end{equation*} 
where $c_0+c_1=9$, $a\in \bC\setminus \Omega$, if\/ $\Omega$ is \underline{unbounded}.
\end{theorem}

\setcounter{theorem}{2}
\begin{theorem}[{on meromorphic functions}] Let\/  $\Omega$ be a subdomain of\/ $\bC$.
Let\/ {\rm (see \eqref{r:f})} $$
f =\frac{p}{q}\in \Mer (\Omega), \; p,q \in \Hol (\Omega), \; p(0)=q(0)=1.$$

Let's choose a function $M$ as in the Theorem\/ {\rm 3 (on zeros)}.  
Suppose that
$$
\sup\limits_{D\in \Exh \Omega} T_{_D}(f) < +\infty.
$$
Then there exist functions  
$p_0,q_0\in \Hol(\Omega; M)$  representing $f=p_0/q_0$.
\end{theorem}

\begin{remark}
If\/ $\Omega$ possesses the Green's function, then we can remove  $\boxed{\sup\limits_{D\in \Exh \Omega}}$ everywhere in the Theorems\/ {\rm 3} and replace $D$ with $\Omega$.
\end{remark}

{\bf General results}. If $M\in \SH (\Omega)$  with the Riesz measure $\nu_M:=\frac1{2\pi}\Delta M \geq 0$, then there is always a global Riesz representation (decomposition)
\begin{equation}\label{con:MR}
	M(z)=\int_{\Omega}  k (\zeta , z) \d \nu_M ( \zeta )+H(z), \; z\in \Omega, 
\end{equation}
where $H\in \Harm (\Omega)$, 
\begin{equation}
k(\zeta, z)=\log |\zeta -z|+h_M(\zeta , z)
\end{equation}
 is a special subharmonic kernel with harmonic component $h(\zeta, z)$ of $z \in \Omega$ for each $\zeta \in \Omega$.
Let $Q \colon \Omega \to \R $ be an upper  
semicontinuous function such that
\begin{equation}	\label{con:MQ}
	\int_\Omega	 \bigl( k(\zeta , 0)-k(\zeta , z )\bigr)^+  \d {\nu}_M (\zeta )
	\leq Q(z) 
\end{equation}
for almost all $z\in\Omega$ in the Lebesgue measure on $\Omega$.  

Denote by $\mathcal U_0^d (\Omega )$ the class of all connected unions $D \ni 0$ of finitely many open disks from $\Omega$  whose complement has no one-point connected components.
\begin{MTheorem}[{see [2, Theorems 1--5]}] Let $M \in \SH (\Omega) \cap C(\Omega)$ with the Riesz measure $\nu_M$.

\begin{enumerate}
	\item[{\bf [Z]}] If $\Lambda=\{ \lambda_k\} \subset \Omega$ is a zero (sub)set for $\Hol(\Omega; M)$, then 
\begin{equation}\label{con:supM}
\sup_{0\in D\Subset \Omega}
\left( \sum_k g_{_D}(\lambda_k, 0)
 -\int g_{_D}(\zeta, 0) \d \nu_M (\zeta)\right)<+\infty .	
\end{equation}
 Conversely, if
 we have  \eqref{con:supM} where domains $D$ run only the class $\mathcal U_0^d (\Omega )$, then $\Lambda$ is zero \underline{sub}sequence (non-uniqueness sequence) for $\Hol (\Omega ; \widehat M )$ where $\widehat M (z):=$
\begin{equation*}
\inf_{0<t<\dist(z, \partial \Omega)}
\Bigl(\frac1{2\pi}\int_0^{2\pi}M\bigl(z+te^{i\theta}
	\bigr)\d \theta +\log \bigl(1+1/t\bigr)\Bigr)+9\log^+|z|,
\end{equation*}
and zero sequence for $\Hol (\Omega ; M +Q)$.
 
Thus, every zero \underline{sub}sequence for $\Hol (\Omega ; M )$ is zero sequence for $\Hol (\Omega ; M +Q)$.
\item[{\bf [M]}] Let $f=g/q$ be a meromorpic function and $g, q \in \Hol (\Omega; M)$. Then 
\begin{equation}\label{con:fM}
\sup_{0\in D\Subset \Omega}
\left( \int_{\Omega} \log \max \{|g|, |q|\}(z)\d \omega_{_D}(0, z) \right. \\ \left. -
\int M(z) \d \omega_{_D}(0, z)\right)<+\infty .	
\end{equation}
 
 Conversely, if, under the assumptions of \eqref{con:MR} and  \eqref{con:MQ}, we have  \eqref{con:fM} where domains $D$ run only the class $\mathcal U_0^d (\Omega )$, then there are functions $g, q \in \Hol (\Omega; \widehat M)$ and $g_0, q_0 \in \Hol (\Omega; M+Q)$ such that $f=g/q=g_0/q_0$ and $g_0, q_0$
have not common zeros. 

Thus, if $f=g/q$ with $g, q \in \Hol (\Omega; M)$, then  there exist functions\\ $g_0, q_0 \in \Hol (\Omega; M+Q)$ without common zeros such that $f=g_0/q_0$.
  
\end{enumerate}
\end{MTheorem}
{\bf Nevanlinna's theorems wit nonradial and nonpositive weight.}
In [4, Theorem 1], we investigate also an slowly counterpart 
Nevanlinna's theorems.

Let $M\colon \D \to \R$, and let $M\in \SH (\D)$ with the Riesz measure $\nu_M$.
Given  $z=re^{i\theta}$, $0\leq  r<1$, $\theta \in \R$, and $a >0$ we consider a polar 
rectangle
\begin{equation}\label{df:sqr}
\hspace{-4mm}	
\boxminus (z;a):=\{ \, \zeta=te^{i\psi}\colon \bigl(r-a\sqrt{1-r^2}\bigr)^+\leq t<1, \; |\sin (\psi-\theta)| <a\sqrt{1-r^2} \, \}
\end{equation}
of relative size $a$ and the function
\begin{equation}\label{df:qma}
q_M^{[a]}(z):=\frac1{1-|z|}\int_{{\boxminus}(z;a) } (1-|\zeta|)\d \nu_M(\zeta).
\end{equation}
We set
\begin{equation}\label{df:aver}
	\Av_M^{[\e]}(z):=\frac1{2\pi}\int_0^{2\pi} M\bigl(z+\e (1-|z|)e^{i\theta}\bigr) \d \theta, \quad 0<\e <1.
\end{equation}
\begin{theorem}\label{th:2.1} Let $M$ be a subharmonic function $\D$, $M(0)>-\infty $ and
\begin{equation}
	\sup_{r<1} 
	\int_0^{2\pi} M(re^{i\theta}) \d \theta <+\infty ,
\end{equation}
that is equivalent a Blaschke condition
\begin{equation}\label{c:BL}
	\int_0^1(1-t) \d \nu_M (t)<+\infty.
\end{equation}

For a function $f_\Lambda$ with $\Zero_{f_\Lambda}=\Lambda$,  
\begin{enumerate}
	\item[{\rm (Z)}] if is fulfilled at least one condition
	\begin{subequations}\label{in:gm1o}
\begin{align}
\sup_{D\in \mathcal U_0^d (\D )} \left( \int_{\D} \log \bigl|f_{\Lambda}(z)\bigr| \d \omega_D(0, z)- \int_{\D} M(z)\d \omega_D(0, z)\right)&< {+\infty},
\tag{\ref{in:gm1o}f}\label{in:gm11}
\\
\sup_{D\in \mathcal U_0^d (\D )} \left(\sum_k g_D(\lambda_k , 0) -
\int_{\D \setminus \{0\}}  g_D (\zeta , 0 )\d \nu_M (\zeta )\right)&< +\infty,\tag{\ref{in:gm1o}g}\label{in:gm1og}\\  
\sup_{D\in \mathcal U_0^d (\D )} \left(\sum_k g_D(\lambda_k , 0) -
\int_{\D} M (z)\d \omega_D(0, z)\right)&< +\infty,\tag{\ref{in:gm1o}h}\label{in:gm1oh}
\end{align}
\end{subequations}
then for any $\e \in (0,1)$ and $1<a <2$ the sequence  $\Lambda \subset \D$ is zero sequence for space
\begin{equation}\label{fin:space}
	\Hol \Bigl(\D; \Av_M^{[\e]}+\frac{C_\e}{2-a}\, q_M^{[a]}\Bigr),  
\end{equation}
where a constant  $C_\e$ dependent only on $\e$;
\item[{\rm (U)}] if $\Lambda$ is zero subsequence for $\Hol (\D ;M)$, then $\Lambda$  is zero sequence for \eqref{fin:space};
\item[{\rm (M)}]  if meromorphic function  $f=g/q$ in\/ $\D$  is represented as a ratio of functions\/ $g, q\in \Hol (\D )$, $\max \bigl\{|g(0)|,|q(0)|\bigr\} \neq 0$, and is fulfilled at least  one condition
\begin{subequations}\label{cond:mer} 
\begin{align}
\hskip -20mm 
\sup_{D\in \mathcal U_0^d (\D )} \left(\int_{\D}\right. \log \left.\max \bigl\{ |g(z)|, |q(z)|\bigr\} \d \omega_D(0, z)- \int_{\D} M(z)\d \omega_D(0, z)\right)&<+\infty,\tag{\ref{cond:mer}h}\label{cond:merh}\\
g, q\in  \Hol (\D ;M )&,
\tag{\ref{cond:mer}M}\label{cond:merM}
\end{align}
\end{subequations}  
then there are functions\/  $g_0$ and $q_0$ from class \eqref{fin:space}
without common zero, such that $f=g_0/q_0$ in $\D$.  
\end{enumerate}
  
\end{theorem}

A function $M\colon z=re^{i\theta} \to [-\infty, +\infty]$ is called {\it radial in a sector
\begin{equation}
	\measuredangle (\alpha ,\beta):=\{z=re^{i\theta} \colon 0\leq r<1, \; \alpha <\theta <\beta \}
\end{equation}
from $\D$,\/} 	if for each $0\leq r <1$ the function $M(re^{i\theta})$ independent of $\theta \in \measuredangle (\alpha , \beta)$.

\begin{corollary} If there is a sector $\measuredangle (\alpha' ,\beta')\Supset \measuredangle (\alpha ,\beta)$ such that  
$\alpha' <\alpha <\beta< \beta'$, and the weight $M$ from Theorem\/ {\rm \ref{th:2.1}} is radial and differentiable in $r$, then, for $a$ from small neighborhood of $1$, for all assertions \/ {\rm (Z),(U) и (M)} of this Theorem at points $z\in  \measuredangle (\alpha ,\beta)$
the summand   $\frac{C_{\e}}{2-a}\,q_M^{[a]}(z)$ in \eqref{fin:space} can be change  to a summand
\begin{equation}
	\frac{aC_{\e}}{2 (2-a)}\, \frac{1}{\sqrt{1-|z|}}
\int_{\bigl(|z|-a\sqrt{1-|z|^2}\bigr)^+}^1 (1-t) \d (tM'(t)).
\end{equation}
for small  $a>1$.
\end{corollary}

This work is supported by the RFSB grants No.~09--01--00046-а, No.~08--01--97023--Volga region, and by the grant of President of Russian Federation ``State support of lead\-ing scientific schools'', project 3081.2008.1.

\renewcommand{\refname}{REFERENCES}

\vsk
 Russia, Ufa, Bashkir State University, khabib-bulat@mail.ru

\vs


\end{document}